\documentclass[journal]{IEEEtran}

\IEEEoverridecommandlockouts

\usepackage{cite}
\usepackage[pdftex]{graphicx}
\usepackage{amsmath}
\usepackage{fixltx2e}
\usepackage{pgfplots}
\usetikzlibrary{matrix,positioning}
\usepackage{tikz}
\usepackage{booktabs, threeparttable}
\usepackage{multirow}
\usepackage{array}

\ifCLASSINFOpdf

\else
 
\fi

\begin{document}

\markboth{This paper has been accepted by the 2018 IEEE Power \& Energy Society General Meeting (PESGM).}%
{Shell \MakeLowercase{\textit{et al.}}: Bare Demo of IEEEtran.cls for IEEE Journals}

\title{Finite-Control-Set Model Predictive Control (FCS-MPC) for Islanded Hybrid Microgrids}


\author{
	\IEEEauthorblockN{Zhehan Yi\IEEEauthorrefmark{1}, Abdulrahman J. Babqi\IEEEauthorrefmark{2}, Yishen Wang\IEEEauthorrefmark{1}, Di Shi\IEEEauthorrefmark{1},
		Amir H. Etemadi\IEEEauthorrefmark{2}, Zhiwei Wang\IEEEauthorrefmark{1}, Bibin Huang\IEEEauthorrefmark{3} }
	\IEEEauthorblockA{\IEEEauthorrefmark{1}GEIRI North America, San Jose, CA 95134}\\
	\IEEEauthorblockA{\IEEEauthorrefmark{2}The George Washington University, Washington, DC 20052} \\
	\IEEEauthorblockA{\IEEEauthorrefmark{3}State Grid Energy Research Institute, Beijing, China} \\
     Email: zhehan.yi@geirina.net
    \thanks{This work is funded by the SGCC Science and Technology Program.}		
}

\maketitle

\begin{abstract}
Microgrids consisting of multiple distributed energy resources (DERs) provide a promising solution to integrate renewable energies, e.g., solar photovoltaic (PV) systems. Hybrid AC/DC microgrids leverage the merits of both AC and DC power systems. In this paper, a control strategy for islanded multi-bus hybrid microgrids is proposed based on the Finite-Control-Set Model Predictive Control (FCS-MPC) technologies. The control loops are expedited by predicting the future states and determining the optimal control action before switching signals are sent. The proposed algorithm eliminates the needs of PI, PWM, and droop components, and offers 1) accurate PV maximum power point tracking (MPPT) and battery charging/discharging control, 2) DC and multiple AC bus voltage/frequency regulation, 3) a precise power sharing scheme among DERs without voltage or frequency deviation, and 4) a unified MPC design flow for hybrid microgrids. Multiple case studies are carried out, which verify the satisfactory performance of the proposed method.    

\end{abstract}


\IEEEpeerreviewmaketitle

\section{Introduction}

\begin{figure}
	\centering
	\includegraphics[width=0.9\columnwidth]{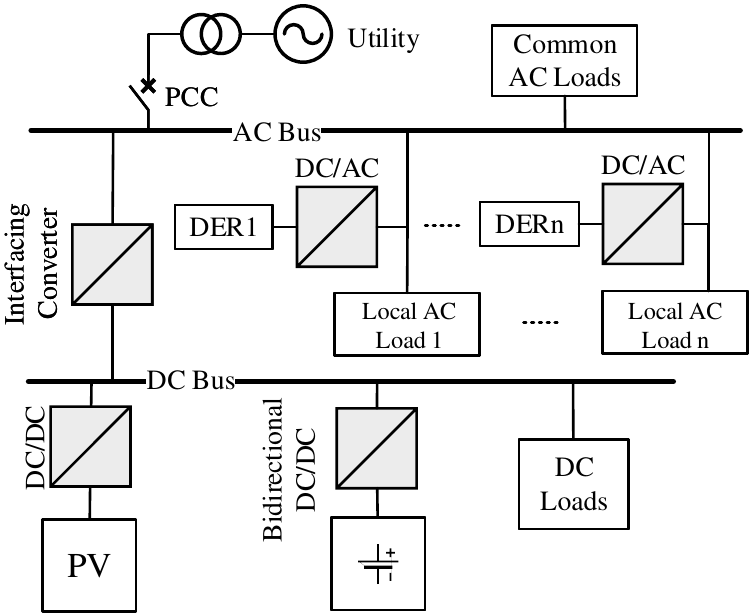}
	\caption{A typical AC/DC hybrid microgrid structure.}\label{structure}
\end{figure}

Deployments of renewable distributed energy resources (DERs) are expected to continuously increase in the near future, due largely to the declining capital investments, political incentives, and attractive natures \cite{TIE,TSG1,pes16}. Microgrid is a promising solution to integrate DERs with advanced control and energy management systems (EMS), which ensures the reliability while simultaneously offers a multitude of benefits to the utility \cite{di,chen}. Introducing a DC network in a conventional AC microgrid, which forms a hybrid microgrid, is advantageous since it allows DC DERs (e.g., PV and battery storage) and DC loads (e.g., LED lighting and DC ventilation) to be interfaced directly with a significantly enhanced efficiency\cite{intro1}. However, the coexistence of both DC and AC buses and DERs requires a sophisticated control system that is capable of managing the power flows effectively and ensuring the stable bus voltages/frequencies. Fig. \ref{structure} illustrates a typical hybrid microgrid, where PV and battery storage are integrated at the DC bus via DC/DC converters and DERs are connected at the AC bus via inverters. Each DER may also include a local bus and load. A bidirectional interfacing converter (IC) is installed between the DC and AC buses. The hybrid microgrid can operate in either grid-connected or islanded mode by switching the circuit breaker at the point of common coupling (PCC). This paper focuses on the control of islanded mode. 

\begin{figure*} 
	\centering
	\includegraphics[width=0.95\textwidth]{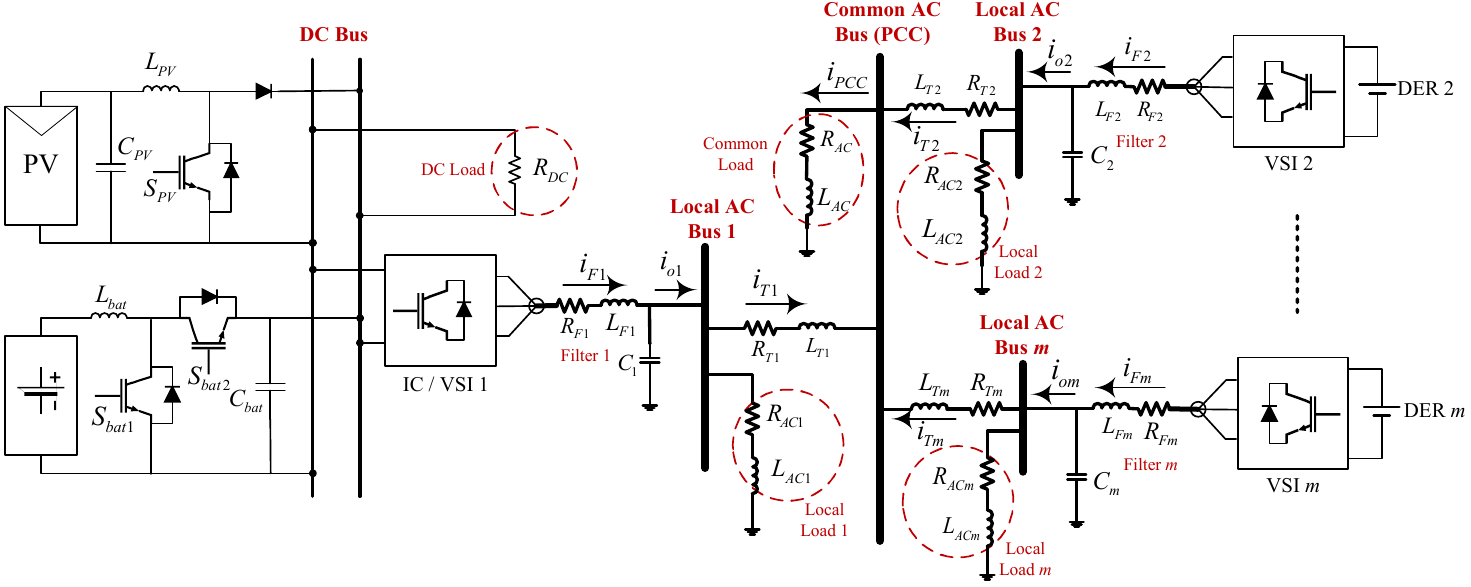}
	\caption{The studied multi-bus AC/DC hybrid microgrid with $m$ VSIs in islanded mode.}\label{simulation_model}
\end{figure*}

Control schemes for conventional microgrids have been extensively studied in the past decade, while those for hybrid microgrids are being actively explored during recent years. To ensure a reliable microgrid, controllers must be able to decouple the buses and the intermittent impacts from renewable DERs, and to control a precise power sharing among multiple DERs. A number of control strategies for hybrid microgrids have been proposed in the literature. Conventional Proportional-Integral (PI)-based controllers with Pulse Width Modulation (PWM) technique have been adopted in a wide range of circumstances, such as in Ref. \cite{pi1,pi2,pi3,pi4,TSG2}. These methods achieve various objectives such as control of multiple interlinking converters, decentralized control without communication links, AC and DC bus regulation, and power management. However, PI technique requires tuning of parameters, which is typically achieved by numerous trial-and-error processes, and needs re-tuning if system changes \cite{abdu}. Moreover, PI controllers should be followed by external PWM modules. Ref. \cite{ann} introduces an adaptive-neural-network (ANN) control scheme for hybrid microgrids, which improves the dynamic performance and reacts adaptively to varying situations by using virtual flux direct power control. Nevertheless, ANN-based methods require massive historic data for the training process and the performance relies largely on the selection of training sets. In terms of power sharing among DERs, droop control and its variations are commonly adopted. A decentralized power management system is proposed in \cite{droop1} to automatically adapt to operating conditions and maintain a power balance, while \cite{pi3} presents a $P_{dc}-V^2_{dc}$ droop control strategy to realize power sharing at the common bus. Similarly, a hierarchical droop-control-based system is presented in \cite{droop2}. Although droop control eliminates the communication links, it suffers from several major drawbacks: 1) voltage or frequency deviation, 2) poor transient performance, and 3) inability of accurate power sharing due to output impedance uncertainties \cite{droop3}. 

As an attempt to address the aforementioned issues, this paper proposes a control scheme for multi-bus hybrid microgrids based on Finite-Control-Set Model Predictive Control (FCS-MPC). The proposed method eliminates the needs of PI, PWM, and droop control, and provides 1) accurate PV maximum power point tracking (MPPT) and battery power control, 2) voltage and frequency regulation with robust setpoint tracking for all buses, 3) a precise power sharing mechanism among the DERs without voltage or frequency deviation, and 4) a unified MPC design flow for hybrid microgrids. The FCS-MPC enables a faster control with satisfactory steady-state and dynamic performance, by predicting the future states of the control objectives and correcting the errors before switching signals are applied to converters/inverters. The rest of paper is organized as follows: Section II elaborates the system modeling and the proposed FCS-MPC for multi-bus hybrid microgrids, case studies based on PSCAD/EMTDC are presented in Section III, and Section IV concludes the paper.

\section{The Proposed FCS-MPC for Hybrid Microgrids}

This section provides a unified process for designing FCS-MPC for multi-bus hybrid microgrids. Fig. \ref{simulation_model} illustrates a typical configuration of an islanded hybrid microgrid with multiple DERs. On the DC side, the PV array and battery storage are connected at the DC bus via DC/DC converters, while on the AC side, each DER is interfaced at the common AC bus (PCC) via a voltage source inverter (VSI). VSI 1 works as an IC to interlink the DC and AC networks. There is a local bus with a load at the output of each VSI and a common bus with a load at PCC. The main philosophy of FCS-MPC is to predict the future behaviors of the system in a predefined time horizon based on the current/past states and possible control actuations. By minimizing a desired cost function, optimal control commands (i.e., switching signals) will be sent, which leads to a minimal error between the objectives and references. Note that the reference for each unit is determined by higher-level EMS, and this paper mainly focuses on primary control. 

The FCS-MPC is designed based on the discrete-time state space of a power electronics stage, which is formulated as: 
\begin{align}\label{state_space}
x(k+1)&=Ax(k)+Bu(k),\\
y(k)&=Cx(k)+Du(k).
\end{align}
\noindent The cost function, Eq. (\ref{cost}), which synthesizes the references, control actuations, and future states of the model, is then minimized subject to certain predefined constraints. 
\begin{equation}
J =  f [x(k),u(k), ..., x(k+N),u(k+N) ] \label{cost}
\end{equation} 
\noindent Optimization process is performed and optimal actuation will be updated as the horizon moves on each sampling time with new samples of measurements \cite{MPC_MPPT, NAPS}. The unified controller design process of each part is elaborated below.

\subsection{PV Controller Design}

\begin{figure} []
	\centering
	\includegraphics[width=0.7\columnwidth]{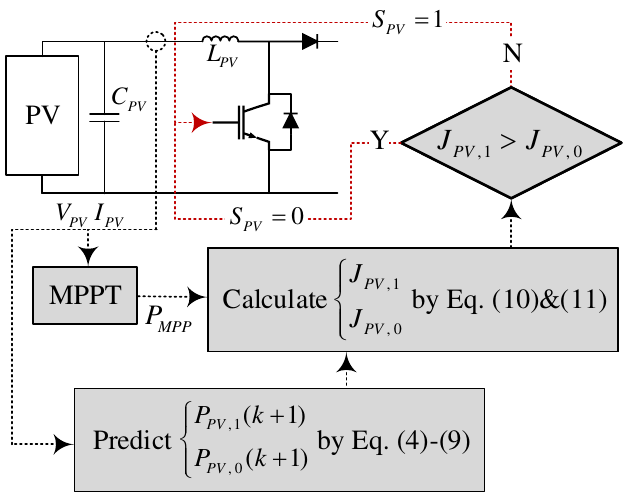}
	\caption{Control algorithm for the PV array.}\label{PV}
\end{figure}

The PV controller aims at extracting the maximum power of the PV array in varying irradiance and temperature. Firstly, a MPPT algorithm (Incremental Conductance \cite{xiaohu}) is employed to determine the real-time maximum power reference ($P_{MPP}$). Based on the state space of DC/DC converter (boost), the next-sample-time value of the PV current ($I_{PV}(k+1)$) and voltage ($V_{PV}(k+1)$) can be derived by Eq. (\ref{Ipv1}) and (\ref{Vpv1}) if the next switch command ($S_{PV}$) is ``ON", and by Eq. (\ref{Ipv0}) and (\ref{Vpv0}) if $S_{PV}$ is ``OFF". $T_{S}$ denotes the sampling period.
\begin{align}
I_{PV,\,1}(k+1)&= I_{PV}(k)+\frac{T_{S}}{L_{PV}}\cdot V_{PV}(k) \label{Ipv1} \\
V_{PV,\,1}(k+1)&= 2V_{PV}(k)-V_{PV}(k-1) \label{Vpv1} \\
I_{PV,\,0}(k+1)&= I_{PV}(k) \label{Ipv0} \\
V_{PV,\,0}(k+1)&= 2V_{PV}(k)-V_{PV}(k-1) \label{Vpv0}
\end{align}
Thereby, the prediction of next PV power is calculated by: 
\begin{align}
P_{PV,\,1}(k+1)&=I_{PV,\,1}(k+1)\cdot V_{PV,\,1}(k+1) \label{Ppv1},\\
P_{PV,\,0}(k+1)&=I_{PV,\,0}(k+1)\cdot V_{PV,\,0}(k+1) \label{Ppv0}.
\end{align}
The cost functions for the PV controller are:
\begin{align}
J_{PV,\,1}&=|P_{PV,\,1}(k+1)-P_{MPP}|, \label{Jpv1} \\
J_{PV,\,0}&=|P_{PV,\,0}(k+1)-P_{MPP}|. \label{Jpv0}
\end{align}
By comparing Eq. (\ref{Jpv1}) and (\ref{Jpv0}), the switching signal that results in a minimal cost will be selected and sent to the converter. This process is presented in Fig. \ref{PV}. Note that the size of PV does not affect the design process.

\subsection{Battery Controller Design}

\begin{figure} []
	\centering
	\includegraphics[width=0.85\columnwidth]{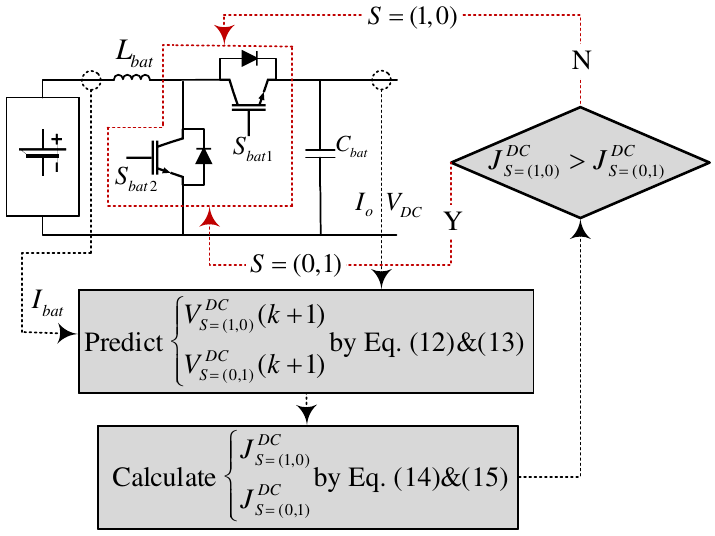}
	\caption{Control algorithm for the battery storage.}\label{bat}
\end{figure}

In the islanded mode of operation, the battery is used to regulate the DC bus while compensating for the power balance between generation and demand. The output voltage of the bidirectional converter, which is also the DC bus voltage, can be predicted by Eq. (\ref{Vdc10}) and (\ref{Vdc01}) when the switching signals $\textbf{S}(S_{bat1}, S_{bat2})$ equal $(1,0)$ and $(0,1)$, respectively.
\begin{align}
V_{\textbf{S}=(1,0)}^{DC}(k+1)&=\frac{T_S}{C_{bat}}\cdot[I_{bat}(k)-I_{o}(k)]+V_{DC}(k) \label{Vdc10} \\
V_{\textbf{S}=(0,1)}^{DC}(k+1)&=-\frac{T_S}{C_{bat}}\cdot I_{o}(k)+V_{DC}(k) \label{Vdc01}
\end{align}
Consequently, the cost functions can be defined by:
\begin{align}
J^{DC}_{\textbf{S}=(1,0)}&=|V^{DC}_{ref}-V_{\textbf{S}=(1,0)}^{DC}(k+1)|, \\
J^{DC}_{\textbf{S}=(0,1)}&=|V^{DC}_{ref}-V_{\textbf{S}=(0,1)}^{DC}(k+1)|.
\end{align}
Depending on the value of each cost function, optimal switching signals that minimize the cost will be delivered to the two switches of the bidirectional converter (Fig. \ref{bat}).

\begin{figure} []
	\centering
	\includegraphics[width=0.7\columnwidth]{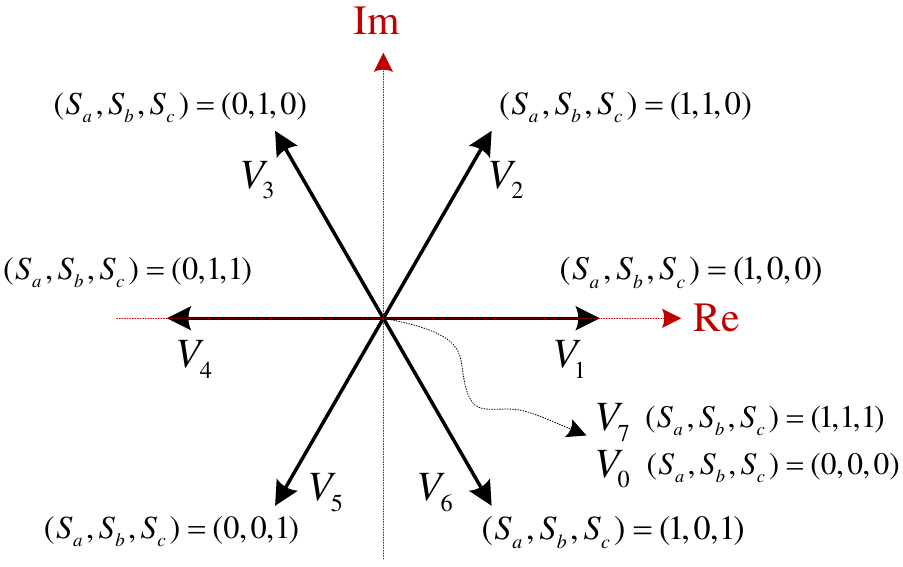}
	\caption{Space Vector Modulation (SVM).}\label{SVM}
\end{figure}

\begin{figure} []
	\centering
	\includegraphics[width=0.75\columnwidth]{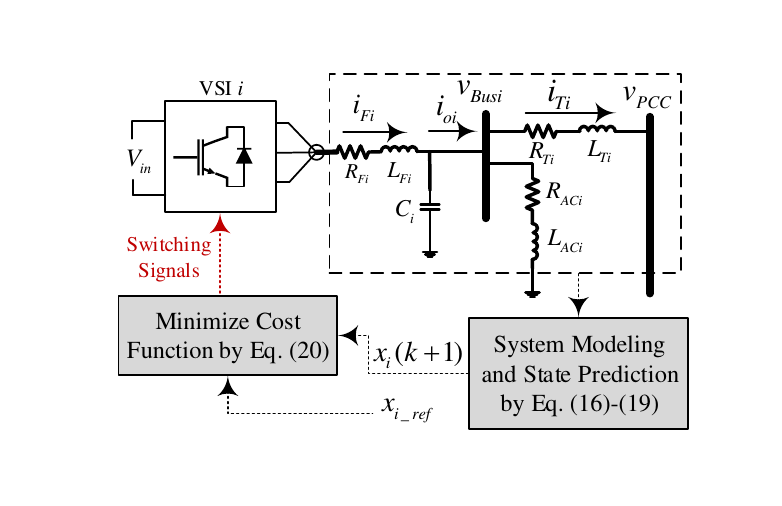}
	\caption{Control algorithm for VSI $i$.}\label{vsi}
\end{figure}

\subsection{VSI Controller Design with Power Sharing Mechanism}

In an islanded multi-DER microgrid system, the voltage and frequency of AC buses are regulated by VSI controllers, which also provide a precise power sharing mechanism among the DERs for the common loads at PCC. The controller design process for each VSI is elaborated below. Each VSI is a conventional three-phase two-level inverter with three legs and two switches in each leg. Therefore, there will be eight possible switching combinations, which yields eight voltage vectors ($V_n$) for the VSI output:
\begin{equation}
V_n=
 \begin{cases}
 \frac{2}{3} V_{in}e^{j(n-1)\frac{\pi}{3}} & n=1,2,...,6\\
 0 & n=0,7
 \end{cases}
\end{equation}

\noindent where $V_{in}$ is the DC input voltage of each VSI. Switching signals for the upper three switches ($S_a$, $S_b$, $S_c$) are given based on the Space Vector Modulation (SVM) technique (Fig. \ref{SVM}), while the lower ones are triggered by ($\bar {S_a},\bar{S_b}, \bar{S_c}$). Therefore, for a microgrid with $m$ VSIs (Fig. \ref{simulation_model}), the state prediction for VSI $i$ is derived as:
\begin{equation}
x_i(k+1)  =A_ix(k)+B_iV_n+C_ii_{oi}(k)
+D_i\sum_{j=1}^{m}i_{Tj}(k)
\end{equation} 

\noindent where
\begin{equation}
\begin{array}{@{} >{{}}r<{{}} @{}r@{\,} *{5}{c} @{\,}l@{}}
x_i = \begin{bmatrix} i_{Fi}\\v_{Busi}\\i_{Ti}  \end{bmatrix},~~~
A_i = \begin{bmatrix} 1-\frac{T_S}{R_{Fi}}  & -\frac{T_S}{L_{Fi}} & 0 \\\frac{T_S}{C_i}  & 1 & 0 \\ 0 & \frac{T_S}{L_{Ti}} & 1-\frac{T_SR_{Ti}}{L_{Ti}} \end{bmatrix},~~~~~~~~~ \\ \\ 
B_i = \begin{bmatrix} \frac{T_S}{L_{Fi}}\\ 0 \\ 0  \end{bmatrix},~~~
C_i = \begin{bmatrix} 0 \\ \frac{-T_S}{C_{i}}\\0 \end{bmatrix},~~~
D_i = \begin{bmatrix} 0 \\0 \\ \frac{-T_S Z_{AC}}{L_{Ti}} \end{bmatrix}.~~~~~~~~~~~~~~~
\end{array}  
\end{equation}

\noindent $v_{Busi}$ denotes the voltage of local AC bus $i$ and $Z_{AC}$ is the load impedance at PCC. The PCC bus voltage becomes:
\begin{equation}
%\begin{split}
v_{PCC}(k+1) = Z_{AC} \sum_{j=1}^{m}i_{Tj}(k+1).
%\end{split}
\end{equation}
Therefore, the optimal space vector that minimizes the following cost function will be selected:
\begin{equation}
\begin{split}
J_{AC} = & \lambda [v_{ref}-v_{PCC}(k+1)]^2 + \\
    & (1-\lambda) \sum_{j=1}^{m-1} [i_{Fj}(k+1)-\beta_j i_{F(j+1)}(k+1)]^2,
\end{split}\label{costfunction}
\end{equation}

\noindent where $v_{ref}$ is the reference for $v_{PCC}$ and $\lambda$ is the weighting factor. The power sharing mechanism is enabled by introducing the second term in $J_{AC}$ with a group of power sharing ratios ($\beta_1,\beta_2,...,\beta_{m-1}$) to control the output current ratio among the $m$ VSIs. For instance, in a system with two VSIs ($m$=2), by setting $\beta_1=1/2$, the output power of VSI 1 will be half of VSI 2. Fig. \ref{vsi} depicts this control process.

\section{Case Studies}

\begin{figure} []
	\centering
	\includegraphics[width=0.9\columnwidth]{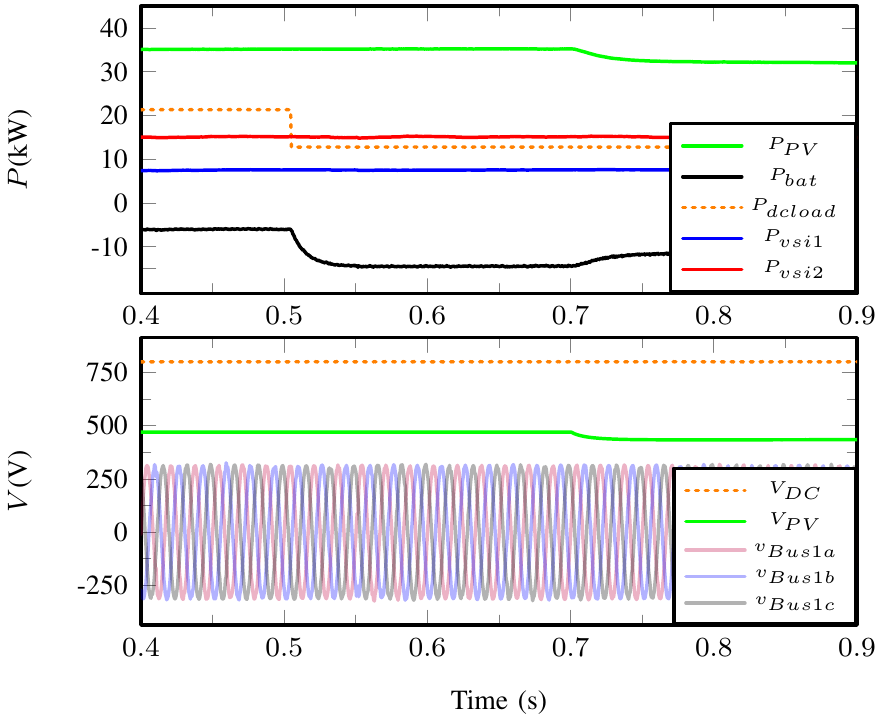}
	\caption{Case 1: powers and voltages during the transitions.}\label{case1}
\end{figure}

To verify the proposed control strategies, a hybrid microgrid with the same configuration as Fig. \ref{simulation_model} with two VSIs ($m$=2) is modeled and the proposed FCS-MPC is implemented by Fortran in the PSCAD/EMTDC software package. Multiple case studies are carried out to validate the performance.

\subsubsection*{Case 1}

This case aims at verifying the proposed scheme for controlling the DC side powers and DC/AC bus voltages. The power sharing ratio is set to $\beta_1=1/2$. As the results presented in Fig. \ref{case1}, at 0.5\,s, the DC load at the DC bus decreases from 21.5\,kW to 13\,kW ($P_{dcload}$). Since the PV array is operating in MPPT mode and the power of VSI 1 is controlled by the MPC, the 8.5\,kW excess power in the DC network will be absorbed by the battery. As a result, the battery charging power ($P_{bat}$) increases from 6\,kW to 14.5\,kW in the negative direction. Note that negative power flow indicates charging the battery \cite{dissertation}. At 0.7\,s, the irradiance of PV array declines. Consequently, PV voltage ($V_{PV}$) drops and PV power ($P_{PV}$) decreases from 35\,kW to 32\,kW. The battery power changes from 14.5\,kW to 11.5\,kW to compensate for the irradiance fluctuation, which ensures a reliable DC load power. The output powers of VSI 1 and 2 remain at 7.5\,kW and 15\,kW ($\beta_1=1/2$), respectively. It is worth noting that transitions of all these changes are finished within 50\,ms, which demonstrates a fast response speed of the proposed scheme. Moreover, the voltages of DC and AC buses are not affected during the transients.

\begin{figure}
	\centering
	\includegraphics[width=0.88\columnwidth]{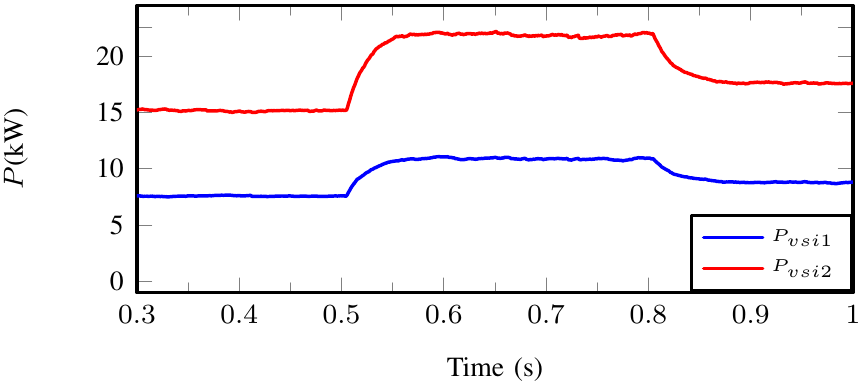}
	\caption{Case 2: output powers of the VSIs with $\beta_1=1/2$.}\label{case2power}
\end{figure}

\begin{figure}
	\centering
	\includegraphics[width=0.9\columnwidth]{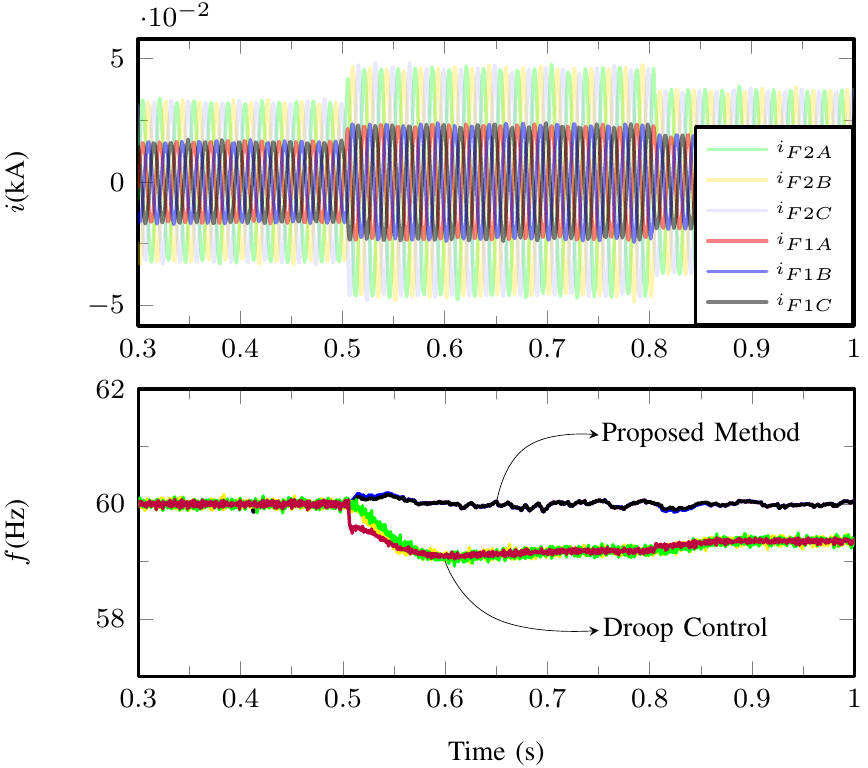}
	\caption{Case 2: 3-phase output currents of the VSIs ($\beta_1=1/2$); comparison of bus frequencies between the proposed method and droop control.}\label{case2current}
\end{figure}

\subsubsection*{Case 2}
The second case study investigates the performance of power sharing mechanism. Since the microgrid is in islanded mode, the PCC common load is shared by VSI 1 and 2, which are initially 7.5\,kW and 15\,kW with $\beta_1=1/2$ (assuming no local loads). As is presented in Fig. \ref{case2power}, at 0.5\,s, the common load at PCC increases by 10.5\,kW. The powers of VSI 1 and 2 increase by 3.5\,kW and 7\,kW, respectively, which results in a 11\,kW output for VSI 1 and a 22\,kW output for VSI 2. Similarly, at 0.8\,s, the common load decreases by 6.6\,kW, resulting in 8.8\,kW for VSI 1 and 17.6\,kW for VSI 2. Note that $P_{vsi2}/P_{vsi1}$ is always 2 as controlled by $\beta_1$. Fig. \ref{case2current} shows the output currents of the two VSIs and the frequencies of three AC buses (Bus 1, Bus 2, and PCC in Fig. \ref{simulation_model}) during these transitions. The current amplitude of VSI 1 is always half of VSI 2. The bus frequencies controlled by the proposed method  60 Hz with negligible fluctuations ($<$0.2\,\%) during the transitions, while extra case studies using droop control show significant bus frequency deviations during load changes.

\subsubsection*{Case 3}

\begin{figure}
	\centering
	\includegraphics[width=0.9\columnwidth]{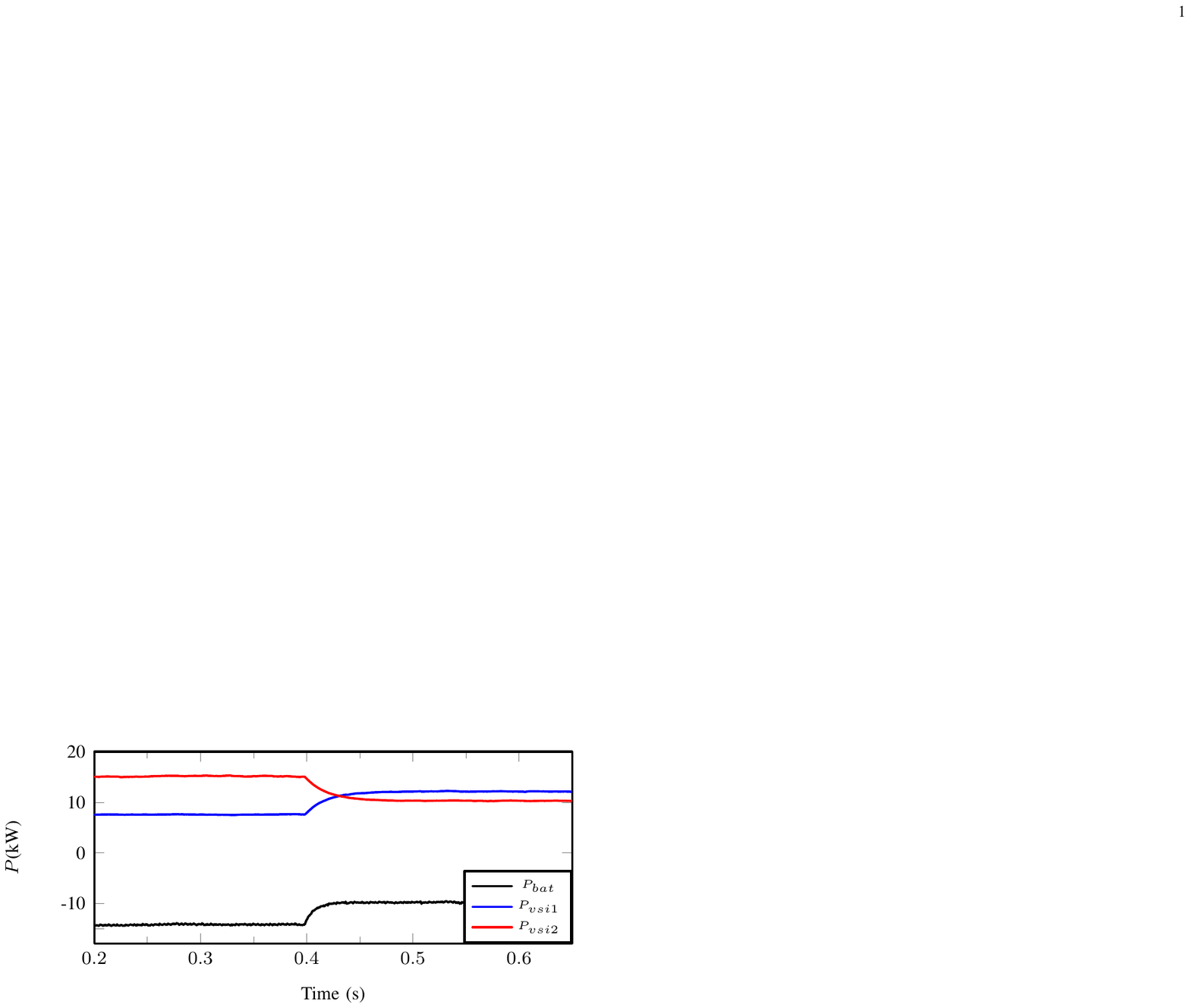}
	\caption{Case 3: power sharing ratio changes from $\beta_1=1/2$ to $\beta_2=8/7$.}\label{case3power}
\end{figure}

This case aims at verifying the flexibility of the proposed method in adjusting power sharing ratio ($\beta_1$). Initially, $\beta_1$ is set to 1/2, which yields 7.5\,kW and 15\,kW for VSI 1 and VSI 2, respectively. At 0.4\,s, $\beta_1$ increases to 8/7. The output of VSI 1 ($P_{vsi1}$ in Fig.\ref{case3power}) increases to 12\,kW whereas the output of VSI 2 ($P_{vsi2}$) decreases to 10.5\,kW, which complies with the change of $\beta_1$ and maintains a reliable power (22.5\,kW) at PCC. Moreover, this verifies how the battery power on the DC side responds to the change.

\section{Conclusions}

This paper introduces a control strategy for multi-bus hybrid microgrids based on FCS-MPC, which eliminates the needs of PI controller, PWM module, and droop control with an improved steady-state and dynamic performance. The proposed scheme predicts the future states of the hybrid microgrid and decides the optimal control actuations before switching signals are sent. It achieves accurate PV MPPT and battery charging/discharging control, DC and AC bus voltage/frequency setpoint tracking, and precise power sharing among DERs at the PCC without voltage or frequency deviation, and offers a unified MPC design approach for hybrid microgrids. Case studies are performed to verify the proposed control strategy.

%\bibliographystyle{IEEEtran}
%\bibliography{pesgm18}

\end{document}